# The relativistic Burgers equation on a de Sitter spacetime. Derivation and finite volume approximation

Baver Okutmustur and Tuba Ceylan

*Abstract*—The inviscid Burgers equation is one of the simplest nonlinear hyperbolic conservation law which provides a variety examples for many topics in nonlinear partial differential equations such as wave propagation, shocks and perturbation, and it can easily be derived by the Euler equations of compressible fluids by imposing zero pressure in the given system. In recent years, several versions of *the relativistic Burgers equations* have been derived on different spacetime geometries by the help of the Lorentz invariance property and the Euler system of relativistic compressible fluid flows with zero pressure on different backgrounds. The relativistic Burgers equation on the Minkowski (flat), Schwarzshild and FLRW spacetime geometries are obtained by LeFloch and his collaborators where the finite volume approximations and numerical calculations of the given models are presented in detail. In this paper, we consider a family member of the FLRW spacetime so-called the de Sitter background, introduce some important features of this spacetime geometry with its metric and derive the relativistic Burgers equation on it. The Euler system of equations on the de Sitter spacetime can be found by a known process by using the Christoffel symbols and tensors for perfect fluids. We applied the usual techniques used for the Schwarzshild and FLRW spacetimes in order to derive the relativistic Burgers equation from the vanishing pressure Euler system on the de Sitter background. We observed that the model admits static solutions. In the final part, we examined several numerical illustrations of the given model through a finite volume approximation on curved spacetimes based on the paper by LeFloch, Amorim and Okutmustur. The effect of the cosmological constant is also numerically analysed in this part. Furthermore, a comparison of the static solution with the Lax Friedrichs scheme is implemented so that the results demonstrate the efficiency and robustness of the finite volume scheme for the derived model

B. O., Middle East Technical University, Department of Mathematics, 06800 Ankara, TURKEY. (phone: 00 90 533 361 62 08, fax: 00 90 312 210 29 72; e-mail: baver@metu.edu.tr).

T. C., Middle East Technical University, Department of Mathematics, 06800 Ankara, TURKEY. (phone: 00 90 505 271 61 56, fax: 00 90 312 210 29 72; e-mail: ceylanntuba@gmail.com).

## I. INTRODUCTION

The classical Burgers equation (with zero viscosity) is one of the simplest nonlinear, hyperbolic partial differential equation model which is important in a variety of applications such as modelling fluid dynamics, turbulence and shocks. It is also the simplest conservation law and formulated by

$$\partial_t u + \partial_x (u^2/2) = 0; \quad u = u(t;x); \, t > 0.$$

In this paper, we derive a relativistic version of the Burgers equation on a de Sitter spacetime and by finite volume approximation, we study the discretization of this equation on the given background. The finite volume method is an important discretization technique for partial differential equations, especially those that arise from conservation laws. This method is the ideal method for computing discontinuous solutions arising in compressible fluid flows and because of the conservation law, they give physically correct weak solutions. In the present paper, we use a finite volume approximation for general balance laws of hyperbolic partial differential equations following the papers [1,3,16,18] and we apply finite volume method to the derived relativistic Burgers equation on a (1+1)-dimensional de Sitter spacetime.

We consider the following class of nonlinear hyperbolic balance laws

$$\text{div}(T(v)) = S(v), \tag{1}$$

on an (n+1)-dimensional spacetime background denoted by M, where the unknown function v is a scalar field, div(.) is the divergence operator, T(v) is the flux vector field and S(v) is the scalar field (source term) represented on M. The spacetime M is assumed to be foliated by hypersurfaces, that is,

$$M = \bigcup_{t \geq 0} H_t, \tag{2}$$

where each slice $H_t$ is an n-dimensional spacelike manifold with initial slice $H_0$. Then the class of nonlinear hyperbolic equations (1)-(2) gives a scalar model on which one can analyze numerical methods. For further details about a description of the curved spacetimes, foliation by hypersurfaces and triangulation of spacelike-timelike elements, we refer the reader to the follow-up works by LeFloch and his collaborators [1,12,13,16].

The inspiration of the relativistic Burgers equation on the de Sitter spacetime is based on the papers [5,13]. We take into account the relativistic Euler equations on a given curved background M by the following general formulation

$$\nabla_\alpha T^{\alpha\beta} = 0, \qquad (3)$$

where $T^{\alpha\beta}$ is the energy-momentum tensor of perfect fluids. The relativistic Burgers equations are derived from (3) by imposing pressure to be zero in this system. The derived equations satisfy the Lorentz invariance property and in the limit case, one obtains the classical (non-relativistic) Burgers equation. In [13], M is firstly taken to be (1+1)-dimensional Minkowski (flat) spacetime and afterwards, it is chosen to be (1+1)-dimensional Schwarzshild spacetime; whereas in [5], the curved spacetime is considered to be the Friedmann-Lemaitre-Robertson-Walker (FLRW) spacetime of (1+1)-dimension. In both of these articles [5,13], the equation (3) is used for deriving the concerning relativistic models for which the numerical tests are taken into account.

In this work, we consider our background to be the de Sitter spacetime which belongs to the family of the FLRW geometry. The metric of both the FLRW and de Sitter backgrounds are solutions of the Einsteins field equations. Our objective is to derive the relativistic Burgers equation on the de Sitter background and examine the numerical results by finite volume approximations for the derived model.

An outline of the article is as follows. We introduce basic features of the de Sitter geometry and its metric in the first part. The derivation of the Euler equations via the Christoffel symbol is in the second part. The next issue will be deriving the particular cases of the relativistic Burgers equations depending on a cosmological constant parameter $\Lambda$. Then we present the finite volume approximations, and compare the classical and the relativistic Burgers equations depending on different values of cosmological constant $\Lambda$. The effects of the cosmological constant to the numerical scheme is also studied in this part. Finally, we terminate the article with numerical experiments by investigating the behaviour of the static solutions and the constructed scheme with different $\Lambda$ values. The results demonstrate that the numerical scheme preserves the static solutions and the scheme is efficient and robust.

II.  DE SITTER BACKGROUND AND ITS LINE ELEMENT

The de Sitter background is a Lorentzian manifold which is a cosmological solution to the Einsteins field equations and it plays an important role in physical cosmology. It shares crucial features with the Minkowski spacetime whereas its physical interpretation is quite complicated. The de Sitter spacetime has a constant curvature and its metric contains a cosmological constant $\Lambda$. The particular case $\Lambda = 0$ in the de Sitter metric gives the Minkowski metric.

The corresponding metric for a (3+1)-dimensional de Sitter spacetime in terms of the proper time t, the corresponding radial r and angular ($\theta$ and $\varphi$) coordinates is formulated by

$$g = -(1-\Lambda r^2)dt^2 + \frac{1}{1-\lambda r^2}dr^2 + r^2(d\theta^2 + \sin^2 d\varphi^2), \qquad (4)$$

which can also be expressed in matrix form as

$$g = g_{ij}dx^i dx^j = (dt\, dr\, d\phi\, d\varphi) D \begin{pmatrix} dt \\ dr \\ d\phi \\ d\varphi \end{pmatrix}$$

where the diagonal matrix D is given by

$$D = \begin{pmatrix} -(1-\Lambda r^2) & 0 & 0 & 0 \\ 0 & \frac{1}{1-\Lambda r^2} & 0 & 0 \\ 0 & 0 & r^2 & 0 \\ 0 & 0 & 0 & r^2 \sin^2\theta \end{pmatrix}.$$

The non-zero **covariant** components and the corresponding **contravariant** components of the diagonal matrix **D** are formulated respectively by

$$g_{00} = -(1-\Lambda r^2), \qquad g_{11} = \frac{1}{1-\Lambda r^2},$$

$$g_{22} = r^2, \qquad g_{33} = r^2 \sin^2\theta,$$

and

$$g^{00} = \frac{1}{\Lambda r^2 - 1}, \qquad g^{11} = 1-\Lambda r^2,$$

$$g^{22} = \frac{1}{r^2}, \qquad g^{33} = \frac{1}{r^2 \sin^2\theta},$$

with

$$g^{ij}g_{kj} = \delta^i_j,$$

where $\delta^i_j$ is the Kronecker's delta function.

### A. Calculation of the Christoffel symbols for a de Sitter background

In order to obtain the Euler equations on a de Sitter background, we start by calculating the Christoffel symbols. The Christoffel symbols are denoted by $\Gamma^\mu_{\alpha\beta}$ and are given by the formula

$$\Gamma^0_{00} = \frac{1}{2} g^{00}(-\partial_0 g_{00} + \partial_0 g_{00} + \partial_0 g_{00}) = 0, \qquad (5)$$

where the terms $\alpha, \beta, \mu, \nu \in \{0,1,2,3\}$. Each term of the Christoffel symbols can easily be calculated by (5). As an example

$$\Gamma^0_{00} = \frac{1}{2} g^{00}(-\partial_0 g_{00} + \partial_0 g_{00} + \partial_0 g_{00}) = 0,$$

$$\Gamma^0_{01} = \frac{1}{2} g^{00}(-\partial_0 g_{01} + \partial_1 g_{00} + \partial_0 g_{10})$$

$$= \frac{1}{2} g^{00}(\partial_1 g_{00}) = \frac{1}{2} \frac{1}{\Lambda r^2 - 1}(\partial_1(-(1-\Lambda r^2)))$$

$$= \frac{1}{2} \frac{1}{\Lambda r^2 - 1}(-(-2\Lambda r)) = \frac{\Lambda r}{\Lambda r^2 - 1}.$$

Repeating this process, the non-zero terms of the Christoffel symbols are obtained as

$$\Gamma^0_{01} = \Gamma^0_{10} = \frac{\Lambda r}{\Lambda r^2 - 1}, \quad \Gamma^1_{11} = \frac{\Lambda r}{1 - \Lambda r^2},$$

$$\Gamma^1_{00} = \Lambda r(\Lambda r^2 - 1), \Gamma^1_{22} = r(\Lambda r^2 - 1),$$

$$\Gamma^1_{33} = r(\Lambda r^2 - 1)\sin^2\theta,$$

$$\Gamma^2_{12} = \Gamma^2_{21} = \Gamma^3_{13} = \Gamma^3_{31} = \frac{1}{r},$$

$$\Gamma^2_{33} = -\sin\theta\cos\theta, \quad \Gamma^3_{23} = \Gamma^3_{32} = \cot\theta.$$

All remaining terms are zero.

### B. Derivation of the relativistic Burgers equation on a de Sitter spacetime

The purpose of this part is to calculate the tensors for perfect fluids on the de Sitter background in order to obtain the Euler equations. We consider our spacetime to be (1+1)-dimensional, that is, the solutions to the Euler equations depend only on the time variable t and radial variable r so that the angular components vanish. Thus we have

$$(u^\alpha) = (u^0(t,r), u^1(t,r), 0, 0), \qquad (6)$$

with

$$u^\alpha u_\alpha = -1,$$

by which we obtain

$$u^\alpha u_\alpha = u^0 u_0 + u^1 u_1,$$
$$= g_{00}(u^0)(u^0) + g_{11}(u^1)(u^1)$$
$$= g_{00}(u^0)^2 + g_{11}(u^1)^2$$
$$= -1$$

It follows by substituting the covariant terms into this relation that

$$-1 = -(1-\Lambda r^2)(u^0)^2 + \frac{1}{\Lambda r^2 - 1}(u^1)^2. \qquad (7)$$

We consider our (3+1)-dimensional coordinates as

$$(x^0, x^1, x^2, x^3) = (ct, r, 0, 0),$$

where c is the light speed. Next, we define the velocity component v by

$$v := \frac{c}{1-\Lambda r^2} \frac{u^1}{u^0}. \qquad (8)$$

Then, by combining the relations (7) and (8), we obtain

$$(u^0)^2 = \frac{c^2}{(1-\Lambda r^2)(c^2 - v^2)},$$

$$(u^1)^2 = \frac{v^2(1-\Lambda r^2)}{(c^2 - v^2)}. \qquad (9)$$

In order to find the tensors, we need the formula of the energy momentum tensor for perfect fluids which is given by

$$T^{\alpha\beta} = (\rho c^2 + p) u^\alpha u^\beta + p g^{\alpha\beta}, \qquad (10)$$

where c is the light speed and $u^\alpha$ is a unit vector described by (7). Then by the help of (9) and (10), the tensors can easily be derived. As an example, the first term can be calculated by

$$T^{00} = (\rho c^2 + p) u^0 u^0 + p g^{00}$$

$$= \frac{c^2(\rho c^2 + p)}{(1-\Lambda r^2)(c^2 - v^2)} + \frac{p}{\Lambda r^2 - 1}$$

$$= \frac{\rho c^4 + p v^2}{(c^2 - v^2)(1-\Lambda r^2)}.$$

Analogously the remaining terms are obtained as follows

$$T^{01} = T^{10} = \frac{cv(\rho c^2 + p)}{(c^2 - v^2)},$$

$$T^{11} = \frac{c^2(1 - \Lambda r^2)(v^2 \rho + p)}{(c^2 - v^2)},$$

$$T^{22} = \frac{p}{r^2},$$

$$T^{33} = \frac{p}{r^2 \sin^2 \theta},$$

$$T^{02} = T^{03} = T^{12} = T^{13} = T^{20} = T^{21} = T^{23}$$
$$= T^{30} = T^{31} = T^{32} = 0.$$

## C. Zero pressure Euler system on a de Sitter background

In the previous section, we derived both the Christoffel symbols and the tensors for the current background. We are now ready to combine these results in order to obtain the Euler system and to impose the pressure to be zero in the concerning equations. To this aim, we use the equation (3), that can be rewritten by

$$\partial_\alpha T^{\alpha\beta} + \Gamma^\alpha_{\alpha\gamma} T^{\gamma\beta} + \Gamma^\beta_{\alpha\gamma} T^{\alpha\gamma} = 0. \quad (11)$$

To begin with $\beta = 0$ in (11), we have

$$\partial_\alpha T^{\alpha 0} + \Gamma^\alpha_{\alpha\gamma} T^{\gamma 0} + \Gamma^0_{\alpha\gamma} T^{\alpha\gamma} = 0.$$

After substituting $\alpha, \gamma \in \{0, 1, 2, 3\}$ in this relation, it follows that

$$\partial_0 T^{00} + \Gamma^0_{00} T^{00} + \Gamma^0_{00} T^{00} + \Gamma^0_{01} T^{10} + \Gamma^0_{01} T^{10}$$
$$+ \Gamma^0_{02} T^{20} + \Gamma^0_{02} T^{02} + \Gamma^0_{03} T^{30} + \Gamma^0_{03} T^{03} + \partial_1 T^{10}$$
$$+ \Gamma^1_{10} T^{00} + \Gamma^0_{10} T^{10} + \Gamma^1_{11} T^{10} + \Gamma^0_{11} T^{11} + \Gamma^1_{12} T^{20}$$
$$+ \Gamma^0_{12} T^{12} + \Gamma^1_{13} T^{30} + \Gamma^0_{13} T^{13} + \partial_2 T^{20} + \Gamma^2_{20} T^{00}$$
$$+ \Gamma^0_{20} T^{20} + \Gamma^2_{21} T^{10} + \Gamma^0_{21} T^{21} + \Gamma^2_{22} T^{20} + \Gamma^0_{22} T^{22}$$
$$+ \Gamma^2_{23} T^{30} + \Gamma^0_{23} T^{23} + \partial_3 T^{30} + \Gamma^3_{30} T^{00} + \Gamma^0_{30} T^{30}$$
$$+ \Gamma^3_{31} T^{10} + \Gamma^0_{31} T^{31} + \Gamma^3_{32} T^{20} + \Gamma^0_{32} T^{32} + \Gamma^3_{33} T^{30}$$
$$+ \Gamma^0_{33} T^{33} = 0.$$

We continue by letting $\beta = 1$ in (11), that is,

$$\partial_\alpha T^{\alpha 1} + \Gamma^\alpha_{\alpha\gamma} T^{\gamma 1} + \Gamma^1_{\alpha\gamma} T^{\alpha\gamma} = 0,$$

and putting $\alpha, \gamma \in \{0, 1, 2, 3\}$ we obtain the following equation

$$\partial_0 T^{01} + \Gamma^0_{00} T^{01} + \Gamma^1_{00} T^{00} + \Gamma^0_{01} T^{11} + \Gamma^1_{01} T^{01}$$
$$+ \Gamma^0_{02} T^{21} + \Gamma^1_{02} T^{02} + \Gamma^0_{03} T^{31} + \Gamma^1_{03} T^{03} + \partial_1 T^{11}$$
$$+ \Gamma^1_{10} T^{01} + \Gamma^1_{10} T^{10} + \Gamma^1_{11} T^{11} + \Gamma^1_{11} T^{11} + \Gamma^1_{12} T^{21}$$
$$+ \Gamma^1_{12} T^{12} + \Gamma^1_{13} T^{31} + \Gamma^1_{13} T^{13} + \partial_2 T^{21} + \Gamma^2_{20} T^{01}$$
$$+ \Gamma^1_{20} T^{20} + \Gamma^2_{21} T^{11} + \Gamma^1_{21} T^{21} + \Gamma^2_{22} T^{21} + \Gamma^1_{22} T^{22}$$
$$+ \Gamma^2_{23} T^{31} + \Gamma^1_{23} T^{23} + \partial_3 T^{31} + \Gamma^3_{30} T^{01} + \Gamma^1_{30} T^{30}$$
$$+ \Gamma^3_{31} T^{11} + \Gamma^1_{31} T^{31} + \Gamma^3_{32} T^{21} + \Gamma^1_{32} T^{32} + \Gamma^3_{33} T^{31}$$
$$+ \Gamma^1_{33} T^{33} = 0.$$

Finally, by substituting the Christoffel symbols and the calculated values of tensors for perfect fluids into the Euler system on (1+1) dimensional de Sitter background, we get a simplified form of the equations as

$$\partial_0 \left( \frac{\rho c^4 + v^2 p}{(1 - \Lambda r^2)(c^2 - v^2)} \right) + \frac{2\Lambda r}{\Lambda r^2 - 1} \left( \frac{cv(\rho c^2 + p)}{(c^2 - v^2)} \right)$$
$$+ \partial_1 \left( \frac{cv(\rho c^2 + p)}{(c^2 - v^2)} \right) + \frac{\Lambda r}{\Lambda r^2 - 1} \left( \frac{cv(\rho c^2 + p)}{(c^2 - v^2)} \right)$$
$$+ \frac{\Lambda r}{1 - \Lambda r^2} \left( \frac{cv(\rho c^2 + p)}{(c^2 - v^2)} \right) + \frac{1}{r} \left( \frac{cv(\rho c^2 + p)}{(c^2 - v^2)} \right)$$
$$+ \frac{1}{r} \left( \frac{cv(\rho c^2 + p)}{(c^2 - v^2)} \right) = 0,$$

$$\partial_0 \left( \frac{cv(\rho c^2 + p)}{(c^2 - v^2)} \right) + \Lambda r(\Lambda r^2 - 1) \left( \frac{\rho c^4 + v^2 p}{(1 - \Lambda r^2)(c^2 - v^2)} \right)$$
$$+ \frac{\Lambda r}{\Lambda r^2 - 1} \left( \frac{c^2(1 - \Lambda r^2)(v^2 \rho + p)}{(c^2 - v^2)} \right)$$
$$+ \partial_1 \left( \frac{c^2(1 - \Lambda r^2)(v^2 \rho + p)}{(c^2 - v^2)} \right)$$
$$+ \frac{2\Lambda r}{1 - \Lambda r^2} \left( \frac{c^2(1 - \Lambda r^2)(v^2 \rho + p)}{(c^2 - v^2)} \right)$$
$$+ \frac{1}{r} \left( \frac{c^2(1 - \Lambda r^2)(v^2 \rho + p)}{(c^2 - v^2)} \right)$$

$$+r(\Lambda r^2 -1)\frac{p}{r^2}+\frac{1}{r}\left(\frac{c^2(1-\Lambda r^2)(v^2\rho+p)}{(c^2-v^2)}\right)$$
$$+r(\Lambda r^2 -1)\sin^2\theta+\frac{p}{r^2\sin^2\theta}=0. \quad (12)$$

Finally we impose p=0 in the above relation (12) in order to get the following two equations so called the **zero pressure Euler system** on (1+1)-dimensional de Sitter background, that is,

$$\partial_0\left(\frac{c}{(1-\Lambda r^2)(c^2-v^2)}\right)+\partial_1\left(\frac{v}{c^2-v^2}\right)$$
$$+\frac{2v}{r(c^2-v^2)}+\frac{2\Lambda rv}{(\Lambda r^2-1)(c^2-v^2)}=0,$$
$$\partial_0\left(\frac{cv}{c^2-v^2}\right)+\partial_1\left(\frac{v^2(1-\Lambda r^2)}{c^2-v^2}\right) \quad (13)$$
$$-\frac{\Lambda rc^2}{c^2-v^2}+\frac{\Lambda rv^2}{c^2-v^2}+\frac{2(1-\Lambda r^2)v^2}{r(c^2-v^2)}=0.$$

## III. THE RELATIVISTIC BURGERS EQUATION ON A DE SITTER BACKGROUND

The objective of this section is to deduce our main equation through the zero pressure Euler equations (13) obtained in the previous section and investigate the static and spatially homogeneous solutions if exist. We combine the first and second equations of (13), with the notation $\partial_0 = \partial_t, \partial_1 = \partial_r$ and derive the following single equation

$$\partial_t v+(1-\Lambda r^2)\partial_r\left(\frac{v^2}{2}\right)+\Lambda r(v^2-c^2)=0, \quad (14)$$

which is the desired **relativistic Burgers equation on a de Sitter background**. We recall that $\Lambda$ is the cosmological constant and c is the speed of light which is a positive parameter.

Depending on various values of $\Lambda$, not only the equation (14), but also the de Sitter metric (4) take different forms. We already stated that, the particular case of the de Sitter metric for $\Lambda=0$ gives the Minkowski metric. In addition, by plugging $\Lambda=0$ in the relativistic Burgers equation (14), one can recover the classical Burgers equation

$$\partial_t v+\partial_r\left(\frac{v^2}{2}\right)=0.$$

The latter observation is a common property shared by the relativistic equations. The limit cases of the relativistic Euler and Burgers equations on the Schwarzschild and FLRW spacetimes have the same feature. For more detail, we refer the reader to the following papers [5,13].

### A. Static and spatially homogeneous solutions

The equation (14) can be written as

$$\partial_t v+\partial_r\left((1-\Lambda r^2)\frac{v^2}{2}\right)=\Lambda r(c^2-2v^2), \quad (15)$$

where the left hand side of the equation is of the conservative form. We firstly investigate the t-independent solutions to the equation (15). To this aim, since the first term of the left hand side $\partial_t v$ vanishes, we consider the following

$$\partial_r\left((1-\Lambda r^2)\frac{v^2}{2}\right)=\Lambda r(c^2-2v^2). \quad (16)$$

In order to solve this equation with respect to r, we use the following change of variable

$$X:=1-\Lambda r^2,$$
$$Y:=c^2-v^2. \quad (17)$$

It follows that
$$Y=KX,$$

where $K\in(0,c)$ is a constant parameter. Combining this result with (17), we get

$$c^2-v^2=K(1-\Lambda r^2),$$

by which, we arrive at a description of the static solutions given by the following formula

$$v_{static}=\pm\sqrt{c^2-K(1-\Lambda r^2)}. \quad (18)$$

On the other hand, it can be observed as a remark that, there is no need for searching spatially homogeneous (r independent) solutions since the equation

$$\partial_t v=\Lambda r(c^2-2v^2)$$

has always r dependency. Therefore, we conclude that our model admits only static solutions.

## IV. FORMULATION OF FINITE VOLUME APPROXIMATION

We consider the hyperbolic balance law (1) in (1+1)-dimensional spacetime M. Thus, the equation (1) turns to be

$$\partial_t T^0(t,r) + \partial_r T^1(t,r) = S(t,r), \qquad (19)$$

where $T^0, T^1$ are flux fields, S is the source term. Following the article by Amorim, LeFloch and Okutmustur [1], we formulate the finite volume approximation by averaging this balance law on the given background.

We intend to write a finite volume scheme on local coordinates for our (1+1)-dimensional background. Supposing that the spacetime is prescribed in coordinates (t,r), we consider the finite volume method over each grid cells

$$[t_n, t_{n+1}] \times [r_{j-1/2}, r_{j+1/2}]$$

with the constant time length $\Delta t$ and the equally spaced cells centered at $r_j$ which are defined by

$$t_n = n\Delta t,$$
$$r_{j+1/2} = r_{j-1/2} + \Delta r,$$
$$r_{j-1/2} = j\Delta r,$$
$$r_j = (j+1/2)\Delta r.$$

In order to obtain the finite volume scheme, we integrate the equation (19) over each grid cell $[t_n, t_{n+1}] \times [r_{j-1/2}, r_{j+1/2}]$ which yields

$$\int_{r_{j-1/2}}^{r_{j+1/2}} (T^0(t_{n+1}, r) - T^0(t_n, r)) dr$$
$$+ \int_{t_n}^{t_{n+1}} (T^1(t, r_{j+1/2}) - T^1(t, r_{j-1/2})) dt$$
$$= \int_{[t_n, t_{n+1}] \times [r_{j-1/2}, r_{j+1/2}]} S(t,r) dt dr.$$

After rearranging the terms it follows that

$$\int_{r_{j-1/2}}^{r_{j+1/2}} T^0(t_{n+1}, r) dr = \int_{r_{j-1/2}}^{r_{j+1/2}} T^0(t_n, r) dr$$
$$- \int_{t_n}^{t_{n+1}} (T^1(t, r_{j+1/2}) - T^1(t, r_{j-1/2})) dt$$
$$+ \int_{[t_n, t_{n+1}] \times [r_{j-1/2}, r_{j+1/2}]} S(t,r) dt dr.$$

Next we introduce the approximations of the numerical flux functions by

$$\tilde{T}_j^n \approx \frac{1}{\Delta r} \int_{r_{j-1/2}}^{r_{j+1/2}} T^0(t_n, r) dr,$$
$$\tilde{Q}_{j\pm 1/2}^n \approx \frac{1}{\Delta t} \int_{t_n}^{t_{n+1}} T^1(t, r_{j\pm 1/2}) dt,$$
$$\tilde{S}_j^n \approx \frac{1}{\Delta r \Delta t} \int_{[t_n, t_{n+1}] \times [r_{j-1/2}, r_{j+1/2}]} S(t,r) dt dr.$$

Thus the finite volume scheme reads as

$$\tilde{T}_j^{n+1} = \tilde{T}_j^n - \frac{\Delta t}{\Delta r} (\tilde{Q}_{j+1/2}^n - \tilde{Q}_{j-1/2}^n) + \Delta t \tilde{S}_j^n, \qquad (20)$$

where $\tilde{T}_j^n = \tilde{T}(v_j^n)$ and $\tilde{Q}_{j\pm 1/2}^n$ are the approximations of the flux functions and $\tilde{S}_j^n$ is the approximation of the source term defined as above. Since the function $\tilde{T}$ is known to be convex (see [13]), by taking the inverse of the above relation, the scheme (20) is rewritten by

$$v_j^{n+1} = \tilde{T}^{-1}\left(T(v_j^n) - \frac{\Delta t}{\Delta r}(\tilde{Q}_{j+1/2}^n - \tilde{Q}_{j-1/2}^n) + \Delta t \tilde{S}_j^n\right).$$

For further details of triangulations and discretization of a geometric and local coordinates formulations of finite volume approximations on a given curved spacetime, we address the reader to the papers [1,13,16] by LeFloch *et al.*

## V. NUMERICAL EXPERIMENTS

### A. Lax Friedrichs scheme for the Burgers equation on a de Sitter background

In this part numerical experiments are illustrated for the model derived on a de Sitter background. The behaviours of initial single shocks and rarefactions are examined in the applications. We use the Lax Friedrichs scheme for finite volume approximations of the derived equation on a de Sitter spacetime by taking $r \in [0,1]$ with various data of cosmological constant $\Lambda$. A standard CFL condition is assumed to be satisfied for the stability of the method.

After normalization (using c=1), the main equation (14) reads as

$$\partial_t v + (1 - \Lambda r^2) \partial_r \left(\frac{v^2}{2}\right) + \Lambda r(v^2 - 1) = 0, \qquad (21)$$

then rewriting the source term on the right hand side of the equation (21), we get

$$\partial_t v + (1 - \Lambda r^2) \partial_r \left(\frac{v^2}{2}\right) = -\Lambda r(v^2 - 1). \qquad (22)$$

The corresponding finite volume scheme for this model is

$$v_j^{n+1} = \frac{1}{2}(v_{j-1}^n + v_{j+1}^n) - \frac{\Delta t}{2\Delta r}(b_{j+1}^n g_{j+1}^n - b_{j-1}^n g_{j-1}^n) + \Delta t S_j^n, \qquad (23)$$

where

$$S_j^n = -\Lambda r_j^n ((v_j^n)^2 - 1),$$
$$b_j^n = (1 - \Lambda (r_j^n)^2), \quad f(v) = v^2/2,$$
$$g_{j-1}^n = f(v_{j-1}^n), \qquad g_{j+1}^n = f(v_{j+1}^n).$$

## B. Effects of cosmological constant $\Lambda$ to the numerical scheme

We examine the effects of the cosmological constant $\Lambda$ to the behaviour of our scheme (23) which is obtained in the previous subsection. To this aim, a comparison of the classical (non-relativistic) Burgers equation and the relativistic Burgers equation derived on a de Sitter background is taken into account. Notice that the classical Burgers equation is a particular case of our model where $\Lambda = 0$. We investigate the attitudes of an imposing shock and rarefaction depending on diverse values of $\Lambda$ on the scheme.

The results are illustrated in the **Figure 1,2,3** and **4**. For each graph, there are two solution curves; the red one represents $\Lambda = 0$ (the classical case) and the green one represents $\Lambda = +1$, -1 (the relativistic cases). The rarefactions are examined in **Figure** 1 and **3**; whereas, shocks are examined in **Figure 2** and **4**. In these experiments we observe that the solution curves related $\Lambda = +1, -1$ move away from the curves corresponding to $\Lambda = 0$ by the time increases.

More precisely, we can realize from the given tests that, for $\Lambda > 0$, the graph of the model extends in upward direction and it is always above the red curve (**Figure 1** and **Figure 2**). On the other hand, for $\Lambda < 0$ the green curve extend always in downward direction and it is below the red curve (**Figure 3** and **4**).

In addition, we observe also that there is an effect of the cosmological constant $\Lambda$ on the speed of the movement of the solution curves. As long as the absolute value of $\Lambda$ becomes larger, the green curve (corresponding to the relativistic one) goes further away from the red one (corresponding to the classical one) in a faster manner.

**Figure 1:** The numerical solutions given by the Lax-Friedrichs scheme with a rarefaction for $\Lambda = 1$.

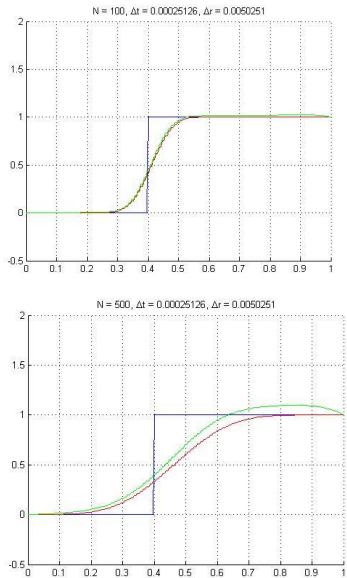

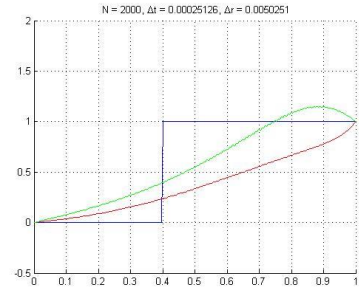

**Figure 2:** The numerical solutions given by the Lax-Friedrichs scheme with a shock for $\Lambda = 1$.

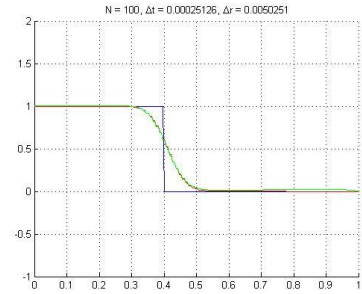

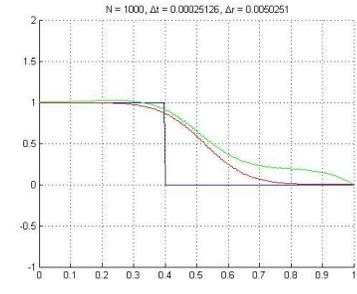

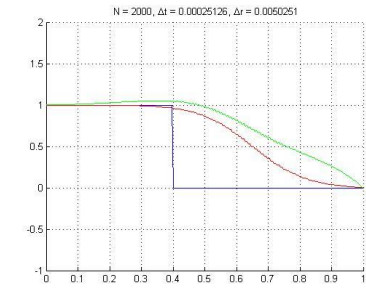

**Figure 3:** The numerical solutions given by the Lax-Friedrichs scheme with a rarefaction for $\Lambda = -1$.

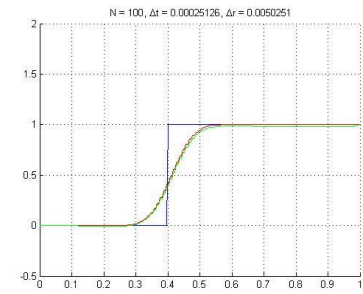

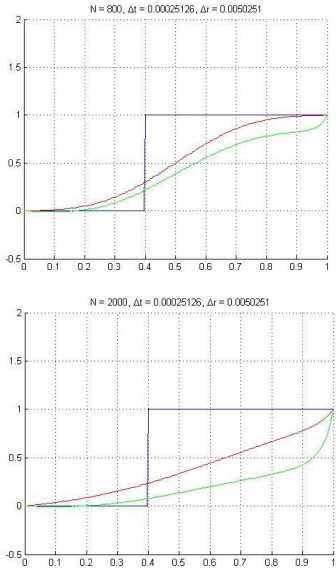

**Figure 4:** The numerical solutions given by the Lax-Friedrichs scheme with a shock for $\Lambda$ = -1.

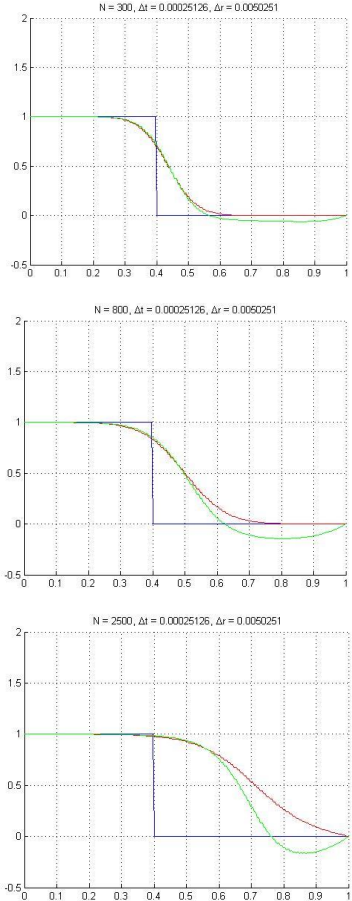

*C. Static solutions for the Burgers equation on a de Sitter background*

We showed previously that our model has static solutions which are described by

$$v_{static} = \pm\sqrt{c^2 - K(1-\Lambda r^2)}.$$

In this part we deal with the numerical tests concerning our scheme and the static solutions. The results of this part are illustrated in **Figure 5, 6** and **7**. We perform the numerical experiments in the domain $r \in [0,1]$ and a usual CFL condition is imposed to be satisfied. A comparison of the static solution curve and the curve resulting from the Lax Friedrichs scheme is under consideration. For each graph, there are two curves; the blue one represents the static solution and the green one represents the constructed scheme.

In **Figure 5**, we start by taking a negative cosmological constant $(\Lambda = -1)$. The initial function is chosen to be the static solution (18) and the constant $K \in [0,1]$ is chosen as *K*=0.5. We observe that the constructed scheme preserves the static solution.

Moreover in **Figure 6**, we take a positive cosmological constant $(\Lambda = 1)$. Analogously, by choosing the static solution as our initial function, we observe that the steady state solution is preserved by the numerical scheme.

Finally in **Figure 7**, we take $\Lambda = 0$. It follows that, the relation (18) yields

$$v_{static} = \pm\sqrt{c^2 - K},$$

which is a constant. For instance, letting (the normalized light speed) c=1 and *K*=0.9, the positive branch of the static solution approximately reads as $v_{static} \approx 0.3162$. By considering this constant as an initial function, the numerical scheme results the same constant curve by which we infer that the constructed scheme again preserves the static solution.

To conclude, after taking into account of positive, negative and zero value cases of cosmological constant, we deduce that the proposed scheme appears to be efficient and robust in the sense that it preserves the static solutions.

**Figure 5:** : Comparison of the static solutions and the constructed schemes with $\Lambda$ = - 1

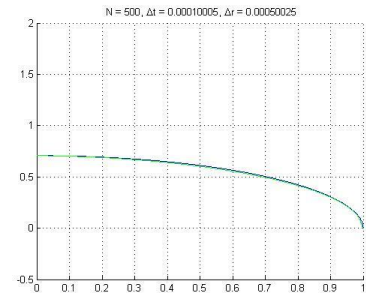

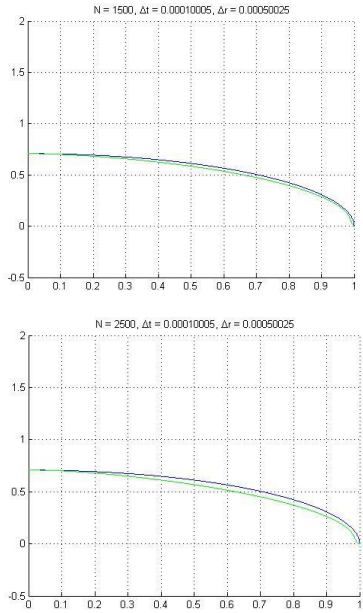

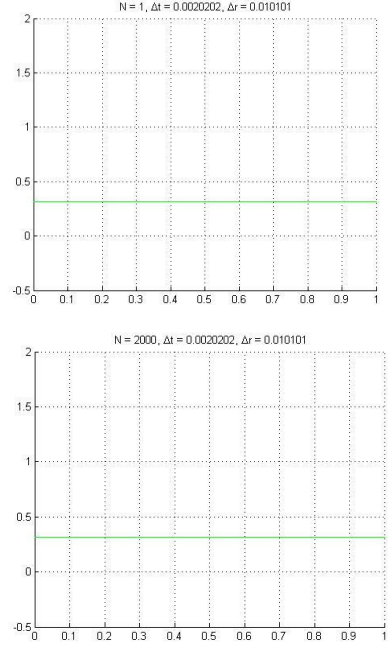

**Figure 7:** Comparison of the static solutions and the constructed schemes with $\Lambda = 0$.

**Figure 6:** Comparison of the static solutions and the constructed schemes with $\Lambda = 1$

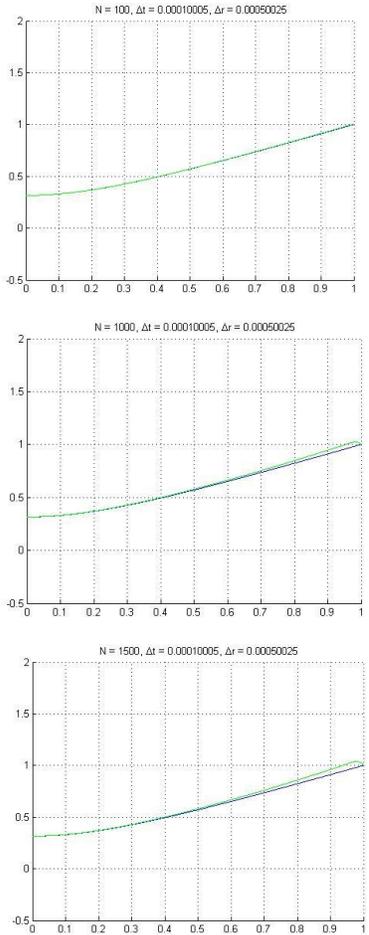

## VI. Conclusion

In this article we derived a new model on a de Sitter background spacetime. We used the technique applied for the Schwarzshild and FLRW spacetimes in order to obtain the relativistic Burgers equation on the de Sitter background. We first used the Cristoffel symbols and tensors for perfect fluids to deduce the concerning Euler system. Then by imposing the pressure to be vanished for this system, we derived the desired relativistic Burgers equation on the de Sitter background.

We observed that the proposed model shares several important features with the relativistic Euler system. Indeed, the unknown v is limited by the light speed parameter and the classical (non-relativistic) model can be derived from the relativistic one similarly as in the Euler system.

Moreover we analysed that for the relativistic Burgers equation on a de Sitter background, the main dependence is in space; in other words, we have static solutions. This led us to examine numerically the behaviour of the constructed scheme in relation with the static solution.

We observed several numerical experiments under two different aspects: the effects of different values of the cosmological constant on the behaviour of the scheme and the relation of the constructed scheme with the static solutions. The first aspect is analyzed by taking into account diverse values of the cosmological constant. We observed that the sign and the magnitude of $\Lambda$ make sense and change not only the direction of the solution curves, but also the speed of these curves. Through the latter aspect, we tested the behaviour of the constructed Lax Friedrichs scheme with respect to the steady state solutions. We demonstrated numerically that the proposed scheme preserves the static solution which guarantees the efficiency and robustness of the constructed finite volume method.

ACKNOWLEDGMENT

The authors would like to thank the reviewers for their valuable comments and suggestions.